\documentclass[11pt]{amsart}
\usepackage{amsmath,amsthm,epsfig,amssymb}
\usepackage{xypic}
\input xy
\title{Counting all cubes  in $\{0,1,...,n\}^3$}
\author{Eugen J. Ionascu and Rodrigo A. Obando}
\curraddr{Department of Mathematics\\ Columbus State University\\4225 University Avenue\\
Columbus, GA 31907\\
} \email{Ionascu\_Eugen@colstate.edu}
\email{Obando\_Rodrigo@colstate.edu}

\date{March $23^{rd}$, 2010}
\begin{document}
\def\sms{\small\scshape}
\baselineskip18pt
\newtheorem{theorem}{\hspace{\parindent}
T{\scriptsize HEOREM}}[section]
\newtheorem{proposition}[theorem]
{\hspace{\parindent }P{\scriptsize ROPOSITION}}
\newtheorem{corollary}[theorem]
{\hspace{\parindent }C{\scriptsize OROLLARY}}
\newtheorem{lemma}[theorem]
{\hspace{\parindent }L{\scriptsize EMMA}}
\newtheorem{definition}[theorem]
{\hspace{\parindent }D{\scriptsize EFINITION}}
\newtheorem{problem}[theorem]
{\hspace{\parindent }P{\scriptsize ROBLEM}}
\newtheorem{conjecture}[theorem]
{\hspace{\parindent }C{\scriptsize ONJECTURE}}
\newtheorem{example}[theorem]
{\hspace{\parindent }E{\scriptsize XAMPLE}}
\newtheorem{remark}[theorem]
{\hspace{\parindent }R{\scriptsize EMARK}}
\renewcommand{\thetheorem}{\arabic{section}.\arabic{theorem}}
\renewcommand{\theenumi}{(\roman{enumi})}
\renewcommand{\labelenumi}{\theenumi}
\newcommand{\Q}{{\mathbb Q}}
\newcommand{\Z}{{\mathbb Z}}
\newcommand{\N}{{\mathbb N}}
\newcommand{\C}{{\mathbb C}}
\newcommand{\R}{{\mathbb R}}
\newcommand{\F}{{\mathbb F}}
\newcommand{\K}{{\mathbb K}}
\newcommand{\D}{{\mathbb D}}
\def\phi{\varphi}
\def\ra{\rightarrow}
\def\sd{\bigtriangledown}
\def\ac{\mathaccent94}
\def\wi{\sim}
\def\wt{\widetilde}
\def\bb#1{{\Bbb#1}}
\def\bs{\backslash}
\def\cal{\mathcal}
\def\ca#1{{\cal#1}}
\def\Bbb#1{\bf#1}
\def\blacksquare{{\ \vrule height7pt width7pt depth0pt}}
\def\bsq{\blacksquare}
\def\proof{\hspace{\parindent}{P{\scriptsize ROOF}}}
\def\pofthe{P{\scriptsize ROOF OF}
T{\scriptsize HEOREM}\  }
\def\pofle{\hspace{\parindent}P{\scriptsize ROOF OF}
L{\scriptsize EMMA}\  }
\def\pofcor{\hspace{\parindent}P{\scriptsize ROOF OF}
C{\scriptsize ROLLARY}\  }
\def\pofpro{\hspace{\parindent}P{\scriptsize ROOF OF}
P{\scriptsize ROPOSITION}\  }
\def\n{\noindent}
\def\wh{\widehat}
\def\eproof{$\hfill\bsq$\par}
\def\ds{\displaystyle}
\def\du{\overset{\text {\bf .}}{\cup}}
\def\Du{\overset{\text {\bf .}}{\bigcup}}
\def\b{$\blacklozenge$}

\def\eqtr{{\cal E}{\cal T}(\Z) }
\def\eproofi{\bsq}

\begin{abstract} In this paper we describe a procedure of calculating the number cubes
that have coordinates in the set $\{0,1,...,n\}$. We adapt the
code that appeared in \cite{ejicregt} developed to calculate the
number of regular tetrahedra with coordinates in the set
$\{0,1,...,n\}$.. The idea is based on the theoretical results
obtained in \cite{ejiam}. We extend then the sequence A098928 in
the Online Encyclopedia of Integer Sequences \cite{OL} to the
first one hundred terms.
\end{abstract} \maketitle
\section{INTRODUCTION}

In this paper we continue and, at the same time revise, some of
the work begun in the sequence of papers \cite{rceji},
\cite{eji}-\cite{ejiam} about equilateral triangles, regular
tetrahedra, and other regular polyhedra, all having integer
coordinates. Very often we will refer to this property by saying
that the various objects are in $\mathbb Z^3$. Strictly speaking
these geometric objects are defined as being more than the set of
their vertices that determines them, but for us here, these are
just the vertices. So, for instance, an equilateral triangle is
going to be a set of three points in $\mathbb Z^3$ for which the
Euclidean distances between every two of these points are the
same. The main purpose of the paper is to take a close look at the
cubes in $\mathbb Z^3$. One can imagine easily such cubes by
taking the faces parallel to the planes of coordinates. However,
it is less obvious that there exist many more other cubes sitting
in space as in Figure~\ref{fig1} (a). As a curiosity, our counting
shows that there are precisely 242,483,634 cubes with vertices in
$\{0,1,...,100\}^3$ and one non-trivial example of these cubes  is
given by the points
$$\{[0, 56, 59], [21, 68, 3], [24, 0, 56], [45, 12, 0], [52, 77, 83],
[73, 89, 27], [76, 21, 80], [97, 33, 24] \}:=\cal C.$$

In \cite{ejiam} we proved the following theorem.
\begin{center}\label{fig1}
$\underset{Figure\ 1 (a): \ Non-trivial\
cube}{\epsfig{file=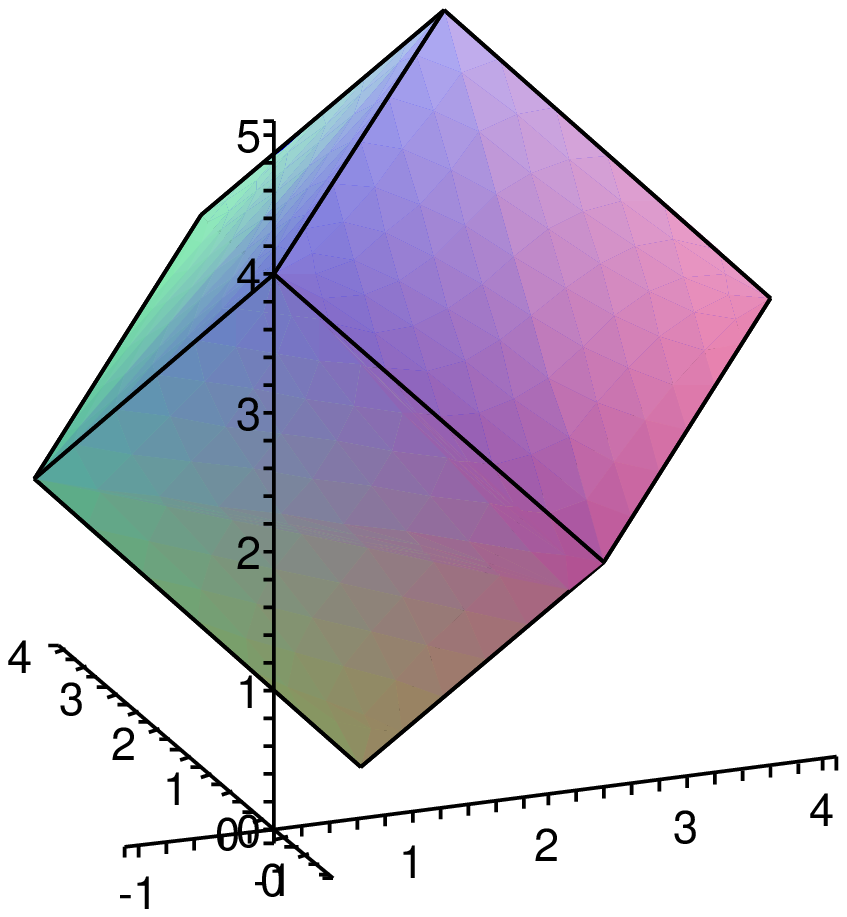,height=3in,width=3in}}\
\underset{Figure \ 1(b): \ Regular\ tetrahedron\ inscribed: \
OABC}{\epsfig{file=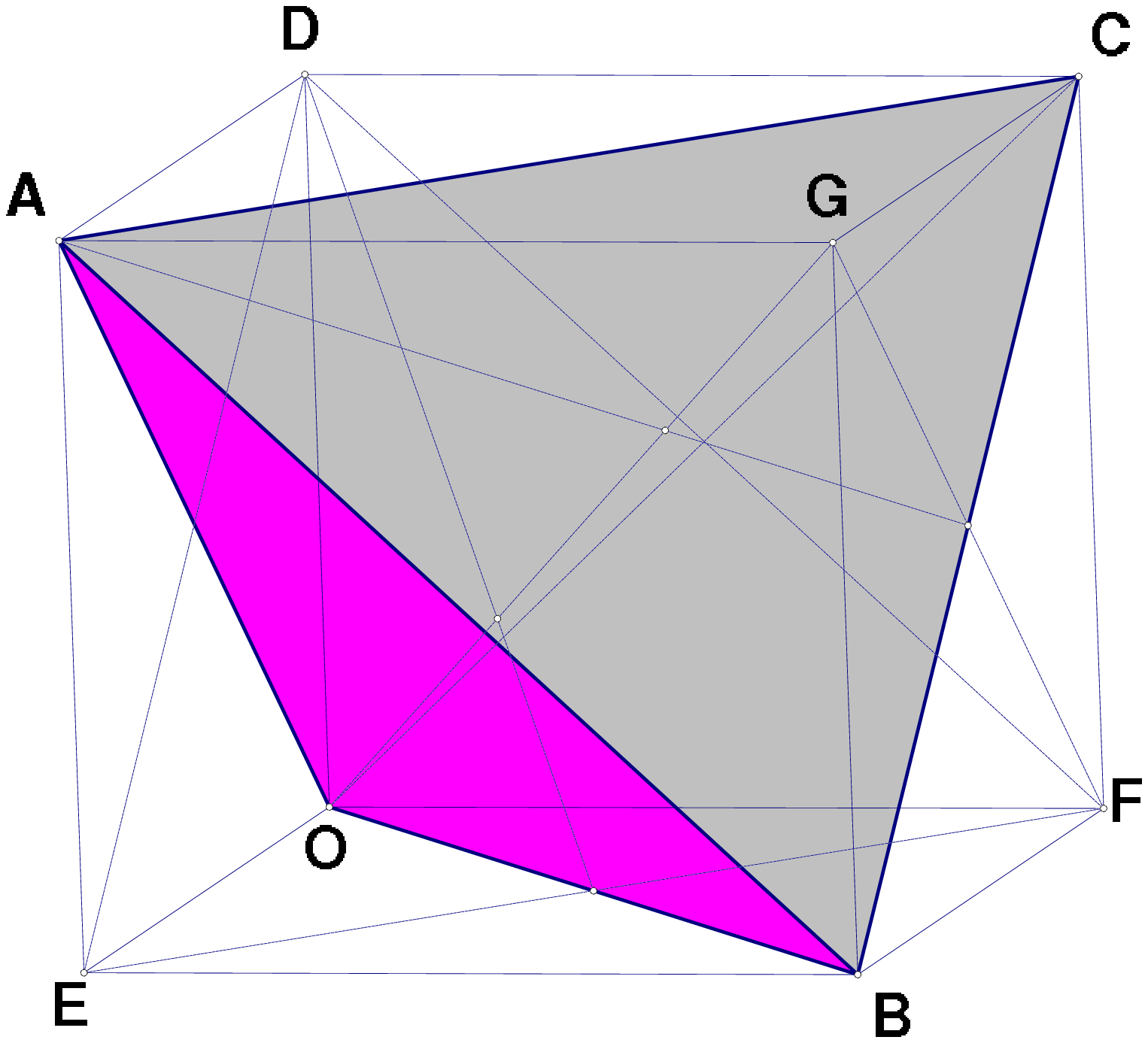,height=2in,width=2in}}$
\end{center}

\begin{theorem}\label{cubes} Every  regular
tetrahedron in $\mathbb Z^3$  can always be completed to a cube in
$\mathbb Z^3$ (See Figure~\ref{fig1} (b)).
\end{theorem}

This theorem implies that there is a one-to-two correspondence
between the cubes and the regular tetrahedrons in $\mathbb Z^3$.
In \cite{ejicregt} we have developed a Maple code to compute the
number regular tetrahedrons in $\{0,1,2,...,n\}^3$. We will
basically use the same idea and introduce some updates based on
important theoretical observations. The problem of finding the
number of cubes in space with coordinates in $\{0,1,2,....,n\}$
has been studied also in \cite{il}. We list here a few more terms
in the sequence A098928.

\vspace{0.1in}

\centerline{ \vspace{0.2in}
\begin{tabular}{|c||c|c|c|c|c|c|c|c|c|c|c|}
  \hline
  n& 1 & 2 & 3 & 4 & 5 & 6 & 7 & 8 & 9 &10&11\\ \hline
 A098928& 1 & 9 & 36 & 100 & 229 & 473 & 910 & 1648 & 2795 & 4469 & 6818  \\
  \hline
\end{tabular}}
\vspace{0.1in}


\centerline{
\begin{tabular}{|c||c|c|c|c|c|c|c|}
  \hline
  n &12&13& 14& 15 & 16 & 17 & 18 \\ \hline
 A098928&10032 & 14315& 19907& 27190 & 36502 & 48233 & 62803 \\
  \hline
\end{tabular}.}
\vspace{0.1in}

It is clear that $A098928 \le A103158$. For $n\ge  4$ we have
actually a strict inequality, $A098928 < A103158$, and this  is
due to the fact that some of the tetrahedrons  inside of the grid
$\{0,1,..,n\}^3$ extend beyond the grid's boundaries  to the
unique cube containing it. In Figure~2 we have included the two
graphs of the sequences A098928 and A103158 up to $n=100$.

\begin{center}
$\underset{Figure\ 2: \ Tetrahedra\ versus \ cubes\
}{\epsfig{file=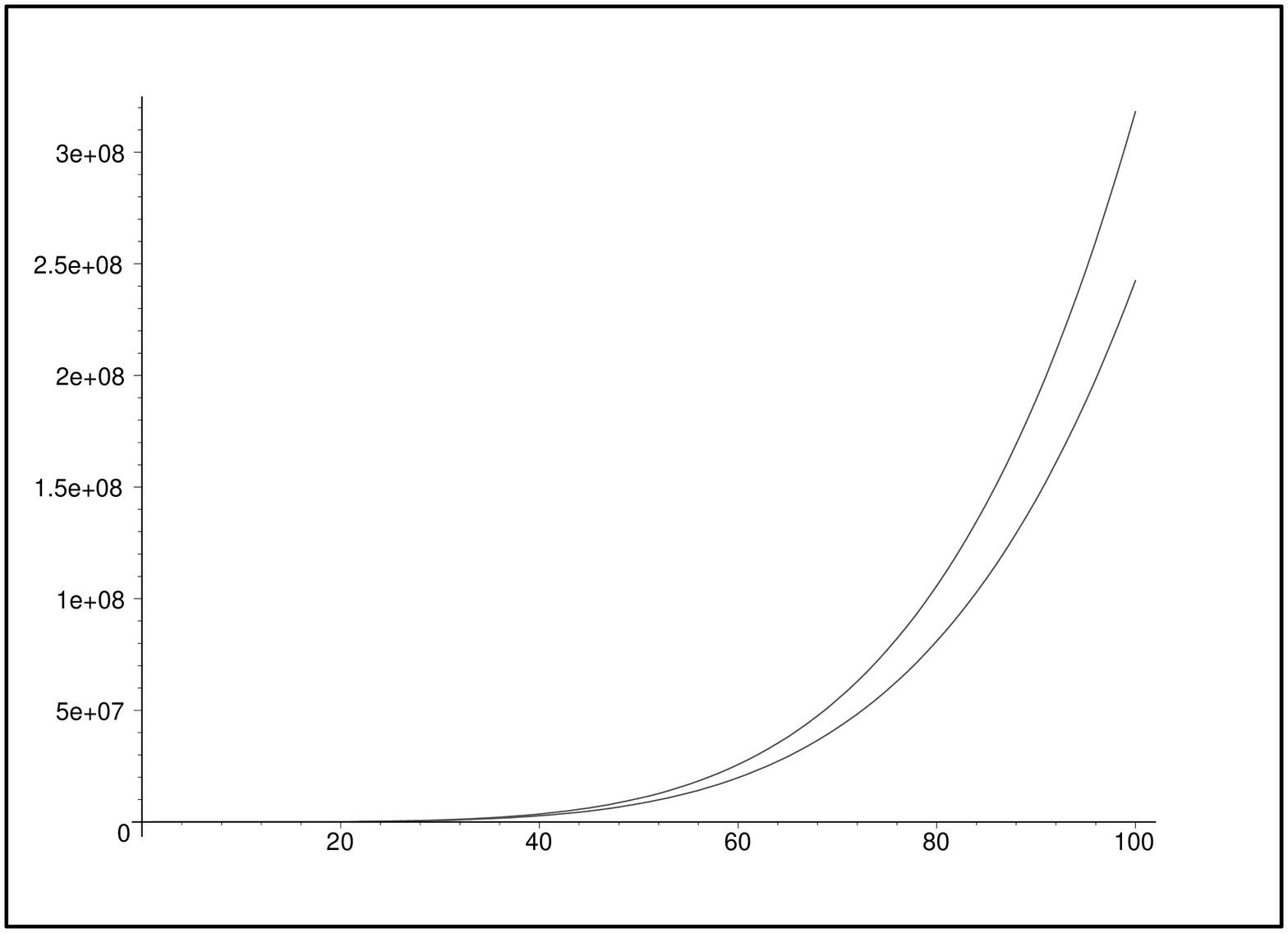,height=3in,width=3.3in,angle=-90}}$
\end{center}

\section{Theoretical background }

Let us review some of the facts that we are using.  Regular
tetrahedra are going to be obtained from equilateral triangles.
Equilateral triangle in $\mathbb Z^3$ are obtained in the
following way. Given an odd integer $d$ there is a very precise
number of solutions for the Diophantine equation (see
\cite{ejicregt}),

\begin{equation}\label{planeeq}
a^2+b^2+c^2=3d^2,\ with\ 0<a\le b\le c\  and \ gcd(a,b,c)=1
\end{equation}

\n which is given by

\begin{equation}\label{numberofrepr}
\pi\epsilon(d)=\frac{\Lambda(d)+24\Gamma_2(d)}{48},
\end{equation}
\n where

\begin{equation}\label{numberofrep2x2plusy2}
\Gamma_2(d)=\begin{cases} 0 \ \text{ if d is divisible by a prime
factor of the form 8s+5 or 8s+7}, \ s\ge 0,\\ \\

1 \ \text{ if d is 3}\\ \\

2^{k}\ \
\begin{cases}\text{ where k is the number of
distinct prime factors of d }   \\ \\
\text {of d of the form 8s+1, or  8s+3}\  (s> 0),\end{cases}
\end{cases}
\end{equation}

\begin{equation}\label{2007hirschhorn}
\Lambda(d):=8d\underset{p|d, p\
prime}{\prod}\left(1-\frac{(\frac{-3}{p})}{p}\right),
\end{equation}

and $(\frac{-3}{p})$ is the Legendre symbol. In particular, these
precise countings show that the equation (\ref{planeeq}) has
solutions for every odd $d\ge 1$.

We remind the reader that, if $p$ is an odd prime then

\begin{equation}\label{legendresymbol}
(\frac{-3}{p})=\begin{cases}0\ \  {\rm if }\  p=3\\ \\
1\ \  {\rm if }\  p\equiv 1\ {\rm or\ } 7 \ {\rm (mod \ 12)} \\ \\
-1\ \  {\rm if }\  p\equiv 5\ {\rm or\ } 11 \ {\rm (mod \ 12)}
 \ \end{cases}.
\end{equation}

As an example, for $d=2011=251(8)+3=167(12)+7$, which is a prime,
will give $\Lambda(d)=16080$ and $\Gamma_2(d)=48$, and so
$\pi\epsilon(2011)=\frac{16080+48}{48}=336$. There is only one
solution, in this case, for which the three values of $a$, $b$,
and $c$ are not all distinct:  $a=139$ and $b=c=2461$.

 For each
solution of (\ref{planeeq}) one can find a two integer parameter
family of equilateral triangle in $\mathbb Z^3$, with vertices in
a plane through any point of integer coordinates (so we can take
simply the origin $O$) and normal $\frac{(a,b,c)}{(d\sqrt{3}}$.
Such an equilateral triangle, say $\triangle OPQ$, can be given
terms of two vectors $\overrightarrow{\zeta}$ and
$\overrightarrow{\eta}$ described by the next formulae.

\begin{equation}\label{vectorid}
\overrightarrow{OP}=m\overrightarrow{\zeta}-n\overrightarrow{\eta},\
\
\overrightarrow{OQ}=n\overrightarrow{\zeta}-(n-m)\overrightarrow{\eta},
\ { \rm with} \ \overrightarrow{\zeta}=(\zeta_1,\zeta_1,\zeta_2),
\overrightarrow{\eta}=(\eta_1,\eta_2,\eta_3),
\end{equation}

\begin{equation}\label{paramtwo}
\begin{array}{l}
\begin{cases}
\ds \zeta_1=-\frac{rac+dbs}{q}, \\ \\
\ds \zeta_2=\frac{das-bcr}{q},\\ \\
\ds \zeta_3=r,
\end{cases}
,\ \
\begin{cases}
\ds \eta_1=-\frac{db(s-3r)+ac(r+s)}{2q},\\ \\
\ds \eta_2=\frac{da(s-3r)-bc(r+s)}{2q},\\ \\
\ds \eta_3=\frac{r+s}{2},
\end{cases}
\end{array}
\end{equation}
where $q=a^2+b^2$ and $(r,s)$ is a suitable solution of
$2q=s^2+3r^2$ that makes all the numbers in {\rm (\ref{paramtwo})}
integers. The sides-lengths of $\triangle OPQ$ are equal to
$d\sqrt{2(m^2-mn+n^2)}$.

One way to give a more precise construction of a good choice of
$(r,s)$ is this to compute the greatest common divisor,
$s+i\sqrt{3}r$, of $A-i\sqrt{3}B$ and $2q$ in the ring $\mathbb
Z[i\sqrt{3}]$, where $A=ac$ and $B=bd$.

Let us observe that $A^2+3B^3=(ac)^2+3(bd)^2=(a^2+b^2)(c^2+b^2)$
which shows that $2q$ divides
$A^2+3B^2=(A+i\sqrt{3}B)(A-i\sqrt{3}B)$. Since $2q=4(4k+1)$ for
some integer $k$, we are thinking of $4$ as
$(1+i\sqrt{3})(1-i\sqrt{3})$, so the prime factors of $2q$ here
are given by $1+i\sqrt{3}$, $1-i\sqrt{3}$ and all the others that
appear which are either primes of the form $6k-1$ or of the form
$6k+1$. The factors of the form $6k-1$ must appear to even power
and those of the form $6k+1$ can be decomposed into prime factors
$u+i\sqrt{3}v$ and $u-i\sqrt{3}v$ (by Euler's Theorem). Each of
these factors can be found either in $A+i\sqrt{3}B$ or
$A-i\sqrt{3}B$. The product of these factors give $s+i\sqrt{3}$.
By construction $A-i\sqrt{3}B=(s+i\sqrt{3}r)(u+i\sqrt{3}v)$ and
$2q=(s+i\sqrt{3}r)(s-i\sqrt{3}r)$ because of properties of
conjugation of complex numbers. This implies that
$(A-i\sqrt{3}B)(s-i\sqrt{3}r))=2q((u+i\sqrt{3}v)$.

Hence we get the relations

$$2q=s^2+3r^2,\ \ As-3Br=2qu\ \ and \ \ Ar+Bs=-2qv.$$

These two relations show that $\zeta_1$ and $\eta_1$ in
(\ref{paramtwo}) are integers. Lagrange's identity shows that

$$\zeta_1^2+\zeta_2^2=\frac{(a^2+b^2)((rc)^2+(ds)^2)}{q^2}=\frac{(rc)^2+(ds)^2}{a^2+b^2},$$

\n which in turn gives

$$\zeta_1^2+\zeta_2^2+\zeta_3^2=\frac{r^2(a^2+b^2)+r^2c^2+d^2s^2}{a^2+b^2}=\frac{d^2(s^2+3r^2)}{q}=2d^2.$$

This implies in particular that  $\zeta_2$ must be an integer and
that $|\overrightarrow{\zeta}|=d\sqrt{2}$.

It is clear that $r$ and $s$ must be either both odd or both even.
This implies that $\eta_3$ is an integer. Using again Lagrange's
identity we get

$$\eta_1^2+\eta_2^2=\frac{(a^2+b^2)[c^2(r+s)^2+d^2(s-3r)^2]}{4q^2}=\frac{c^2(r+s)^2+d^2(s-3r)^2}{4(a^2+b^2)},$$

\n which implies

$$\eta_1^2+\eta_2^2+\eta_2^2=\frac{(r+s)^2(a^2+b^2+c^2)+d^2(s-3r)^2}{4(a^2+b^2)}=\frac{d^2[3(r+s)^2+(s-3r)^2]}{4q}=2d^2.$$

As before, this proves that $\eta_2$ is an integer and
$|\overrightarrow{\eta}|=d\sqrt{2}$.

In order to find the dot product of $\overrightarrow{\zeta}$,
$\overrightarrow{\eta}$ we observe that

$$\zeta_1\eta_1+\zeta_2\eta_2=\frac{(a^2+b^2)[c^2(r^2+rs)+d^2(s^2-3rs)]}{2q^2}$$

\n which implies

$$\overrightarrow{\zeta}\cdot
\overrightarrow{\eta}=\frac{c^2(r^2+rs)+d^2(s^2-3rs)}{2q}+\frac{r^2+rs}{2}=\frac{d^2(3r^2+3rs+s^2-3rs)}{2q}=d^2.$$

Hence the angle between the vectors $\overrightarrow{\zeta}$,
$\overrightarrow{\eta}$ is
$\arccos(\frac{\overrightarrow{\zeta}\cdot \overrightarrow{\eta}}
{|\overrightarrow{\zeta}| | \overrightarrow{\eta}|})=60^{\circ}$.

Using these relations, we can easily calculate
$$|\overrightarrow{OP}|^2=m^2|\overrightarrow{\zeta}|^2-2\overrightarrow{\zeta}\cdot
\overrightarrow{\eta}mn+n^2|\overrightarrow{\eta}|
^2=2d^2(m^2-mn+n^2),\ \ and$$

$$|\overrightarrow{OQ}|^2=n^2|\overrightarrow{\zeta}|^2-2\overrightarrow{\zeta}\cdot
\overrightarrow{\eta}n(n-m)+(n-m)^2|\overrightarrow{\eta}|
^2=2d^2[n^2-n(n-m)+(n-m)^2]=2d^2(m^2-mn+n^2).$$

The dot product of $\overrightarrow{OP}$ and $\overrightarrow{OQ}$
is then equal to

$$\begin{array}{l}
\overrightarrow{OP}\cdot
\overrightarrow{OQ}=mn|\overrightarrow{\zeta}|^2-[m(n-m)+n^2]\overrightarrow{\zeta}\cdot
\overrightarrow{\eta}+n(n-m)|\overrightarrow{\eta}| ^2=\\
\\ d^2(2mn-mn+m^2-n^2+2n^2-2mn)=d^2(m^2-mn+n^2).
\end{array}$$

These relations show that the triangle $\triangle OPQ$ is indeed
equilateral and its side lengths are equal to
$d\sqrt{m^2-mn+n^2}$. One can easily check that

$$a\zeta_1+b\zeta_2+c\zeta_3=a\zeta_1+b\eta_2+c\eta_3=0,$$

\n which implies that $\triangle OPQ$ is indeed contained in the
plane of normal $\overrightarrow{n}=\frac{(a,b,c)}{d\sqrt{3}}$.

One good question here is, whether or not there are other
equilateral triangles with integer coordinates contained in the
same plane.

We have proved this in \cite{rceji}. However, for completion we
will include a new relatively simpler argument and with a more
geometric flavor. If there exists one such triangle, say
$\triangle OAB$, we may assume one of its vertices is at the
origin (by implementing a translation with integer coordinates).
Because the vectors $\overrightarrow{\zeta}$ and
$\overrightarrow{\eta}$ form a basis for the space of vectors
perpendicular to $\overrightarrow{n}$, the equation
$\overrightarrow{OA}=m\overrightarrow{\zeta}-n\overrightarrow{\eta}$
can be solved uniquely for real numbers  $m$ and $n$. Let us
consider then the vectors
$$\overrightarrow{OB'}=n\overrightarrow{\zeta}-(n-m)\overrightarrow{\eta},\ \
\overrightarrow{OB''}=(m-n)\overrightarrow{\zeta}-m\overrightarrow{\eta}.$$
With the same computations as before this will give two
equilateral triangles in the same plane, namely $\triangle OAB'$
and $\triangle OAB''$. Because there are only two equilateral
triangles sharing the side $OA$ in the given plane, we must have
either $B'=B$ or $B''=B$. Without loss of generality, let us
assume that $B'=B$. From the formulae in (\ref{paramtwo}) we get
that

$$mr-\frac{r+s}{2}n=u\in \mathbb Z,\ \ \ and\ \ \
nr-\frac{r+s}{2}(n-m)=\frac{r+s}{2}m+\frac{r-s}{2}n=v\in \mathbb
Z.$$ If we look at these two equations a system of equations in
$m$ and $n$, we get by Cramer's formula, a unique solution which
are going to be (rational numbers) fractions with integer
numerator as the denominators equal to
$$r\frac{r-s}{2}+\frac{r+s}{2}\frac{r+s}{2}=\frac{s^2+3r^2}{4}=\frac{q}{2}=\frac{a^2+b^2}{2}\ge 1.$$
The same calculations as before show that the side length of
$\triangle OAB$ is given by $\ell=d\sqrt{2(m^2-mn+n^2)}$. We can
do this whole construction for $b$ and $c$, or for $a$ and $c$,
instead of $a$ and $b$. We get that

\begin{equation}\label{side}
\ell^2=2d^2(m^2-mn+n^2)=2d^2(m_1^2-m_1n_1+n_1^2)=2d^2(m_2^2-m_2n_2+n_2^2)
\end{equation}

\n with $m_1$, $n_1$, $m_2$, $n_2$ rational numbers with
$\frac{b^2+c^2}{2}$ and $\frac{a^2+c^2}{2}$ at the denominators,
respectively. But the
$gcd(\frac{a^2+b^2}{2},\frac{b^2+c^2}{2},\frac{a^2+c^2}{2})$
cannot be greater than one since $gcd(a,b,c)=1$. Hence, in
(\ref{side}) we cannot have the number
$m^2-mn+n^2=m_1^2-m_1n_1+n_1^2=m_2^2-m_2n_2+n_2^2$ be a fraction
in the reduced form with a denominator greater then one since any
prime dividing it will divide
$gcd(\frac{a^2+b^2}{2},\frac{b^2+c^2}{2},\frac{a^2+c^2}{2})=1$.
So, we proved that the triangle $\triangle OAB$ (or any other
triangle with integer coordinates in the plane containing the
origin of normal $\overrightarrow{n}$) has sides at least
$d\sqrt{2}$.

Now, if $m$ and $n$ are not integers, then $A$ and $B$ fall
strictly inside of the tessellation with equilateral triangles
generated by the two vectors $\overrightarrow{\zeta}$ and
$\overrightarrow{\eta}$ (see Figure 3). Because the tessellation
is invariant to $60^{\circ}$ the position of $B$ inside one of the
equilateral triangles is perfectly the same as the position of $A$
inside of its equilateral triangle. This creates two vectors of
the length and one is the rotation of the other by $60^{\circ}$.
Using translations with integer coordinates the two vectors show
the existence of a triangle with one vertex the origin, $\triangle
OCD$, which is equilateral having integer coordinates and side
lengths strictly less than $d\sqrt{2}$. This contradiction shows
that $A$ and $B$ must be vertices of the tessellation generated by
$\overrightarrow{\zeta}$ and $\overrightarrow{\eta}$, and so the
parametrization (\ref{paramtwo}) is unique.
\begin{center}\label{tessellations}
$\underset{Figure\ 3: \ Two\ distinct \ tessellations\
}{\epsfig{file=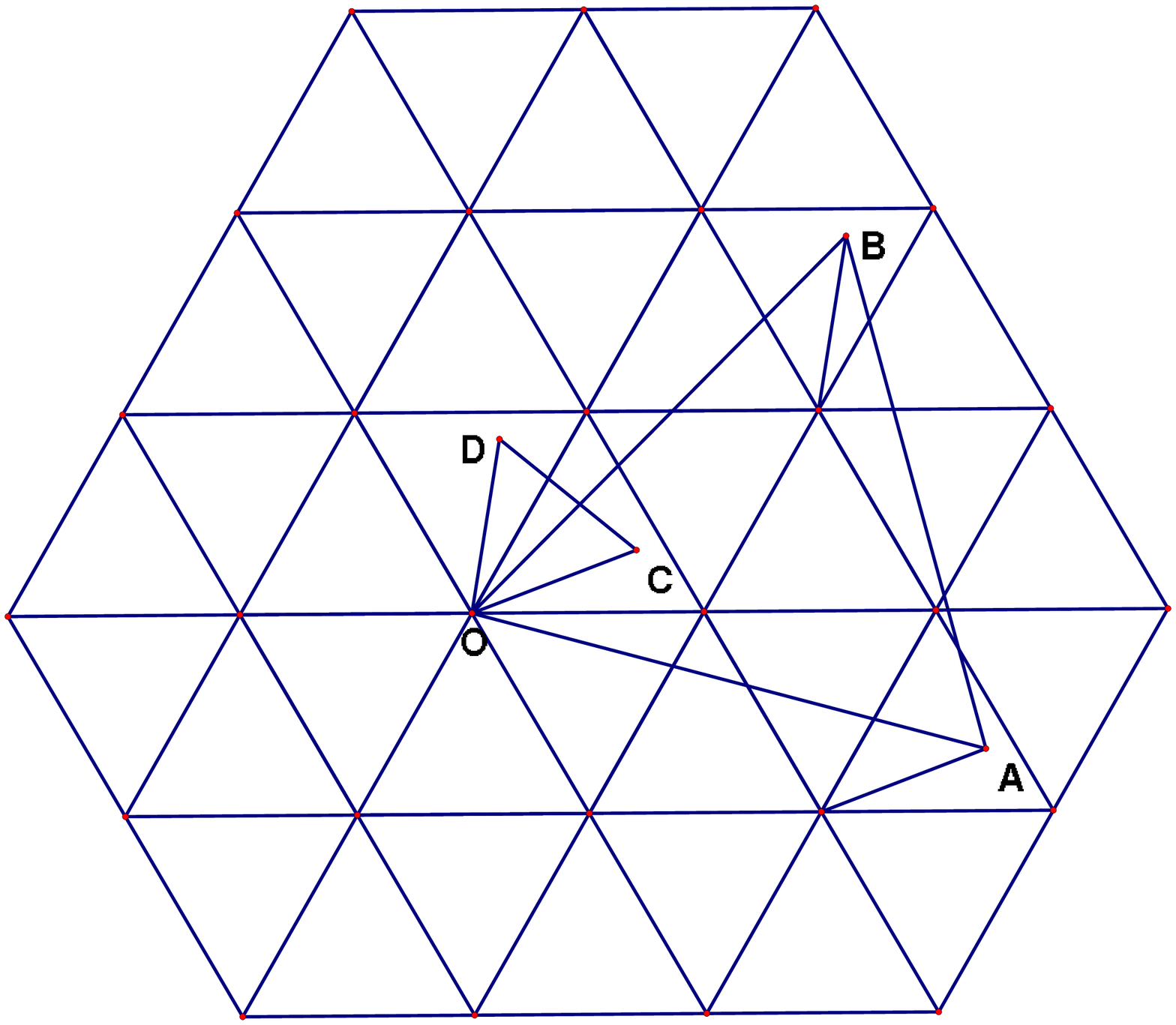,height=3in,width=3.3in}}$
\end{center}

We have given then another proof of Theorem~1 in \cite{ejic}.

We are including here an example that is illustrating this
parametrization. For $d=2011$ we have seen a solution of
(\ref{planeeq}): $a=139$ and $b=c=2461$. If we do the
parametrization with $q=a^2+b^2=(2)(3037921)$ (3037921 is prime).
Since $A=ac=(23)(107)(139)$ and $B=(23)(107)(2011)$, we get

$$A-iB=(23)(107)(1-\sqrt{3}i)(1543-468\sqrt{3}i)$$

and

$$2q=(1-\sqrt{3}i)(1543-468\sqrt{3}i)(1+\sqrt{3}i)(1543+468\sqrt{3}i)$$

which gives
$s+r\sqrt{3}i=(1-\sqrt{3}i)(1543-468\sqrt{3}i)=139-2011\sqrt{3}i$

This gives

$$\begin{cases}
\zeta_1=0\\
\zeta_2=2011\\
\zeta_3=-2011
\end{cases} \ and\  \ \
\begin{cases}
\eta_1=-2461\\
\zeta_2=1075\\
\zeta_3=-936.
\end{cases}$$

One can check that, in fact, in the case $b=c$ we may always take

\begin{equation}\label{bequalc}
\begin{cases}
\zeta_1=0\\
\zeta_2=d\\
\zeta_3=-d
\end{cases} \ and\  \ \
\begin{cases}
\eta_1=-b\\
\eta_2=\frac{a+d}{2}\\
\eta_3=\frac{a-d}{2}.
\end{cases}
\end{equation}

Let us summarize all these facts that we have shown so far.

\begin{theorem}\label{main1} For every solution of the
Diophantine equation $a^2+b^2+c^2=3d^2$ there exists essentially
only one parametrization for the equilateral triangles with
vertices in the plane of equation $ax+by+cz=0$. Every equilateral
triangle in $\mathbb Z^3$ is given by such a parametrization.
Moreover, up to a translation, the vertices of such a triangle are
given by (\ref{vectorid}) and (\ref{paramtwo}), where $r$  and $s$
can be computed by finding the greatest common divisor,
$s+i\sqrt{3}r$, of $A-i\sqrt{3}B$ and $2q$ in the ring $\mathbb
Z[i\sqrt{3}]$, where $A=ac$ and $B=bd$. In the case $b=c$, the
formulae (\ref{paramtwo}) simplify to (\ref{bequalc}).
\end{theorem}

In \cite{ejirt}, we have shown that the only equilateral
triangles, in $\mathbb Z^3$, which can be completed to a regular
tetrahedron in $\mathbb Z^3$, are the ones (given as in
(\ref{vectorid}) and (\ref{paramtwo})) for which $m^2-mn+n^2=k^2$
for some $k\in \mathbb Z$. More precisely, if $k$ is divisible by
$3$ then one can accomplish this on either side of the plane
containing the triangle and if $k$ is not divisible by $3$ then
this can be done on only one side. By the way, this is saying in
particular that, there are a lot more equilateral triangles  than
regular tetrahedrons in $\mathbb Z^3$.

The coordinates for the fourth vertex, assuming the equilateral
triangle's vertices are as in (\ref{vectorid}) and
(\ref{paramtwo}), are given by

\begin{equation}\label{fourthpoint}
\begin{array}{c}
 \displaystyle \left(
\frac{\begin{array}{c}
        (2\zeta_1-\eta_1)m \\
        -(\zeta_1+\eta_1)n \\
       \pm 2ak
      \end{array}
}{3},\frac{\begin{array}{c}
        (2\zeta_2-\eta_2)m \\
        -(\zeta_2+\eta_2)n \\
       \pm 2bk
      \end{array}}{3},
 \displaystyle \frac{\begin{array}{c}
        (2\zeta_3-\eta_3)m \\
        -(\zeta_3+\eta_3)n \\
       \pm 2ck
      \end{array}}{3} \right).
\end{array}
\end{equation}

As we already mentioned in Theorem~\ref{cubes}, as long as the
coordinates in (\ref{fourthpoint}) are integers then the
tetrahedron can be completed to a cube in $\mathbb Z^3$. We are
using this formula mostly for $k=1$ (let us choose $m=1$ and
$n=0$) although there is a need for the general case for big
values of $d$ because, as pointed out in \cite{ejicregt}, there
are irreducible regular tetrahedra which cannot be constructed
from a face as above, by simply taking $k=1$. However for small
$\ell$ ($\ell <5187=3(7)(13)(19)$) one can find a face of a given
regular tetrahedron of sides equal to $\ell \sqrt{2}$ which has
the corresponding $k$ as in (\ref{fourthpoint}) equal to $1$.

\begin{center}
\label{alleightregtetrahedra}
$\underset{Figure\ 4: \ Eight\ tetrahedra \ and \ essentially\
one\ cube\ }{\epsfig{file=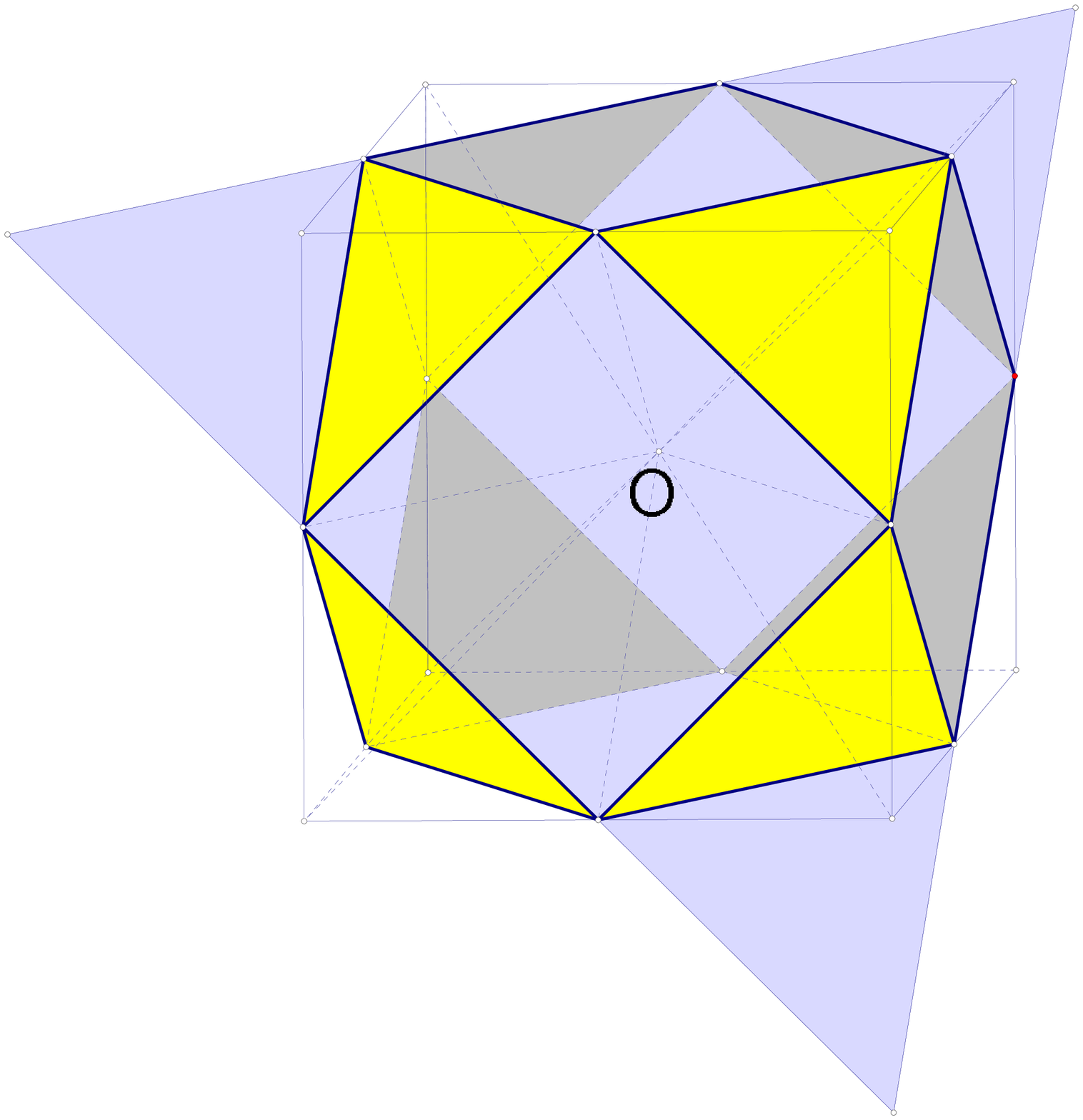,height=3in,width=3.3in}}$
\end{center}

If one takes all the possible values for $m$ and $n$ such that
$m^2-mn+n^2=1$ there are six regular tetrahedra generated this
way, from a plane (colored blue in Figure~4), three on one side
and the other three on the other side,  but if one looks at the
Figure~4, he/she might observe that in fact there are eight
regular tetrahedra all generating essentially the same cube (up to
translations of integer coordinates). As result, our code is going
to reflect this property and we will only use one of the choices
for the values for $m$ and $n$. Summarizing, there are in general
four planes containing the center of a given cube in $\mathbb
Z^3$, corresponding to normals given by the directions of the four
big diagonals in the cube which may generate the the cube as
before, some may have a value of $k>1$. For this reason, one needs
to check for repetitions when writing the code. For this purpose,
our approach is to generate an exhaustive  list, $\cal L$, of
cubes in $\mathbb N_0^3$ ($\mathbb N_0=\mathbb N\cup \{0\}$) which
are irreducible (cannot be integer dilated to a smaller cube in
$\mathbb Z^3$). One other property of each cube in $\cal L$ is
that it cannot be translated in the negative direction along any
of the axes of coordinates and remain in $\mathbb N_0^3$. However,
the cubes in $\cal L$ are not uniquely defined this way, because
of the possible symmetries involved here. These symmetries are, in
general, 48 in number and form a group which can be identified
with the symmetry group of a regular octahedron (see
\cite{ejiam}).

\section{The minimal List and other considerations}

\begin{center}\label{cubestogether}
$\underset{Figure\ 5:\  First\ cubes\ in\ \cal
L}{\epsfig{file=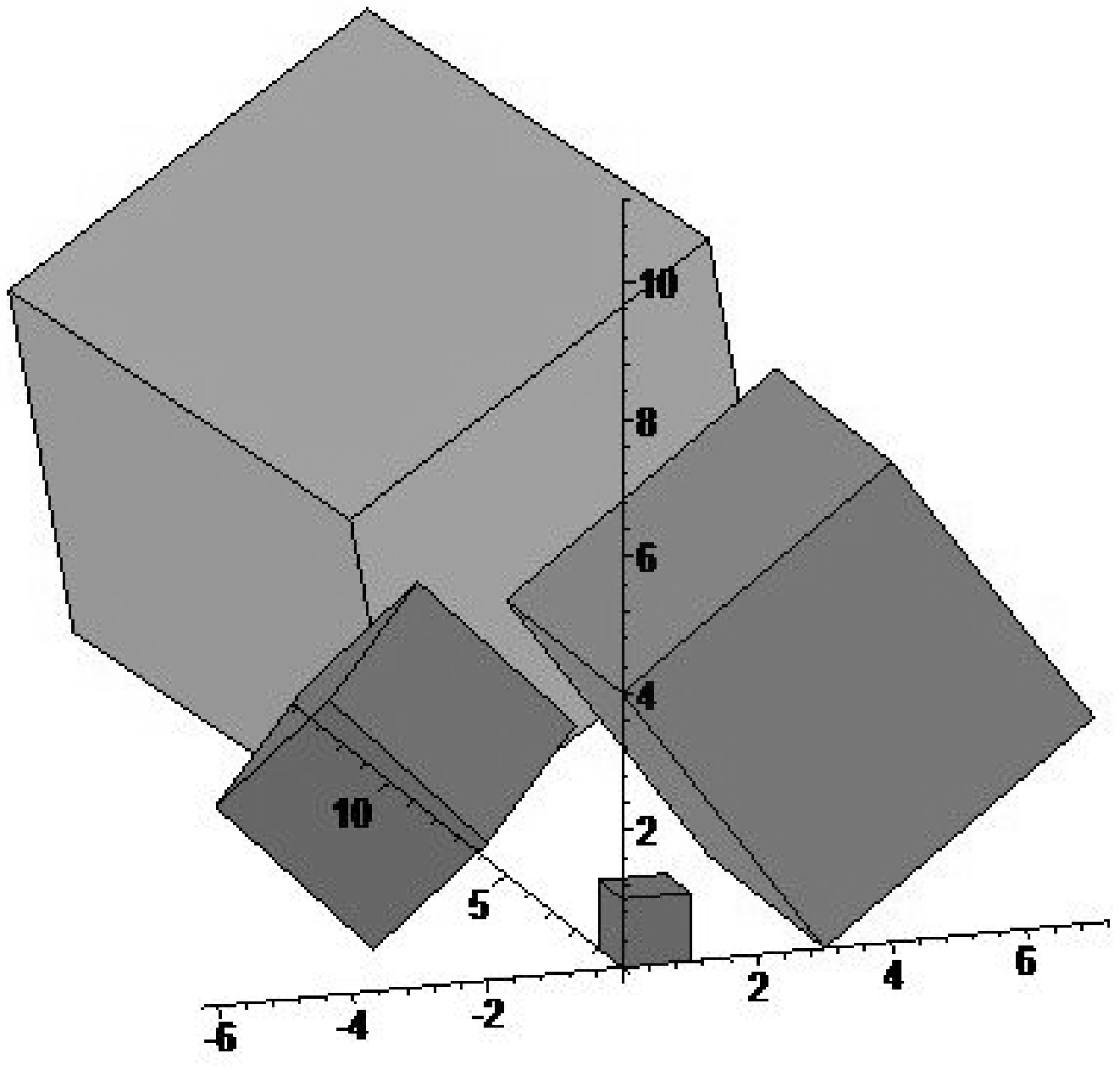,height=3in,width=3in}}$
\end{center}

A dozen cubes (listed in non-decreasing order of their
side-lengths) in the list $\cal L$ that we found using our code,
are included in the table below. The first column represents the
side-lengths, the second column gives the dimension of the
smallest cube $C_m:=[0,m]^3$ containing the one in column three.

\vspace{0.1in}

{ \n  \centerline{Table 1}\vspace{0.1in}
\par \n \tiny
\begin{tabular}{|c|c|c|c|c|}
  \hline   \hline
 n & m & A cube & k-values & invariants \\
   \hline   \hline
  1 & 1 & {\tiny [0, 0, 0], [0, 0, 1], [0, 1, 0], [0, 1, 1], [1, 0, 0], [1,
0, 1], [1, 1, 0], [1, 1, 1]}  & 1 &  [1, 1, 0, 0] \\
  \hline
  3 & 5 & {\tiny [0, 3, 2], [1, 1, 4], [2, 2, 0], [2, 5, 3], [3, 0, 2], [3, 3, 5],
[4, 4, 1], [5, 2, 3]}  & 1,3 & [4, 4, 0, 0] \\
  \hline
  5 & 7 & {\tiny [0, 0, 4], [0, 5, 4], [3, 0, 0], [3, 5, 0], [4, 0, 7], [4, 5, 7],
[7, 0, 3], [7, 5, 3]}  & 1 & [12, 18, 4, 0]
\\  \hline
7 & 11 &{\tiny [0, 6, 8], [2, 9, 2], [3, 0, 6], [5, 3, 0], [6, 8,
11], [8, 11, 5], [9, 2, 9], [11, 5, 3]} & 1,7 &[8, 8, 0, 0] \\
\hline 9& 15& {\tiny [0, 5, 5], [4, 4, 13], [4, 13, 4], [7, 1, 1],
[8, 12, 12], [11, 0, 9], [11, 9, 0], [15, 8, 8]}& 1,3 &{\tiny [24,
108, 48, 16]}\\ \hline 11& 19& {\tiny [0, 11, 13], [2, 2, 7], [6,
17, 6], [8, 8, 0], [9, 9, 19], [11, 0, 13], [15, 15,
12], [17, 6, 6]}& 1 &{\tiny [24, 108, 48, 16]} \\
\hline 13& 19& [0, 12, 15], [3, 16, 3], [4, 0, 12], [7, 4, 0],
[12, 15, 19], [15, 19, 7], [16, 3, 16], [19, 7, 4]& 1, 13 & [8, 8,
0, 0]\\
\hline 13& 17& [0, 0, 12], [0, 13, 12], [5, 0, 0], [5, 13, 0],
[12, 0, 17], [12, 13, 17], [17, 0, 5], [17, 13, 5]& 1&  [12, 30,
8, 0]\\ \hline 15& 25& [0, 5, 10], [2, 19, 15], [10, 0, 20], [11,
7, 0], [12, 14, 25], [13, 21, 5], [21, 2, 10],[23, 16, 15]&  1, 3&
[48, 360, 176, 64]\\ \hline 17& 29& [0, 20, 9], [1, 8, 21], [12,
12, 0], [12, 29, 17], [13, 0, 12], [13, 17, 29], [24, 21, 8], [25,
9, 20]& 1& [24, 60, 16, 0]\\ \hline 17& 23& [0, 0, 15], [0, 17,
15], [8, 0, 0], [8, 17, 0], [15, 0, 23], [15, 17, 23], [23, 0, 8],
[23, 17, 8]& 1 & [12, 42, 12, 0]\\ \hline 19& 31& [0, 16, 10], [6,
6, 25], [10, 31, 16], [15, 10, 0], [16, 21, 31], [21, 0, 15], [25,
25, 6], [31, 15, 21]& 1, 19& [8, 8, 0, 0]\\ \hline
 \end{tabular}}
\vspace{0.1in}

\n  In the column four we list the values of $k$ which can be used
in the construction described in Section 2 to generate the cube in
column three. The list of invariants are as follows.  First, is
the number of cubes in the orbit obtained by applying the group of
48 transformations determined by the orthogonal matrices with
coefficients 0 and 1. Let us denote this number by $\alpha_0$.
Notice that this is a divisor of 48 as expected (Lagrange's
Theorem). We expect that in general such a cube will have no
special symmetry and so, more often we will get $\alpha_0=48$. The
second number in the invariants list is the number of cubes in the
{\it generalized} orbit, obtained by the previous orbit together
with all its integer translations along the axes of coordinates
that remain in $C_m$, a number that we are going to denote by
$\alpha$. The third number in the list, $\beta$, is the
cardinality of the intersection between this former orbit and its
translation along $(0,0,1)$. Finally, the last number, $\gamma$,
is defined by the cardinality of the generalized orbit with its
translation along $(0,-1,1)$. It turns out that these last  three
numbers are enough to determine the number of cubes that one can
fit by translating the given cube in all possible ways within a
bigger cube of size $k\ge m$. This fact has been essentially
proved in Theorem 2.2 in \cite{ejic}. The formula that gives this
number is

\begin{equation}\label{oldformula}
(k-m+1)^3\alpha-3(k-m)(k-m+1)^2\beta+3(k-m+1)(k-m)^2\gamma.
\end{equation}
One of the observations that we will make, about Table 1, is that
this set of invariants is not complete, since we see that the same
numbers appear for various irreducible cubes. The most surprising
are those cubes in rows six and seven. A good problem here is to
determine the exact number of such cubes, which go into a certain
side-length $n$, in terms of $n$. We see that the first $n$ for
which we have two such cubes is $n=13$. Let us also observe that
in column four we see a 1 in there for each cube. As we mentioned
earlier, this is not always the case.

Also, each cube in Table 1 with side-lengths $n$, gives rise to an
orthogonal matrix with rational coefficients having denominators
in $\frac{1}{n}\mathbb Z$ (obtained by taking the normalized
vectors along the sides of  the cube that form an orthogonal
basis). In \cite{ejiam} we computed a few of them which are
included here next:
 \vspace{0.1in}

{\small
$$T_3:= \frac{1}{3}
\left[ \begin{array}{rrr}
          1 & -2 & 2 \\
          2 & -1 & -2 \\
          -2 & -2 & -1 \\
        \end{array}
      \right], T_5:= \frac{1}{5}
\left[ \begin{array}{rrr}
          4 & 0 & 3 \\
          3 & 0 & -4 \\
          0 & -5 & 0 \\
        \end{array}
      \right],
      T_7:= \frac{1}{7}
\left[ \begin{array}{rrr}
          -2 & -3 & 6 \\
          3 & -6 & -2 \\
          -6 & -2 & -3 \\
        \end{array}
      \right],
      T_9:= \frac{1}{9}
\left[ \begin{array}{rrr}
          -7 & -4 & 2 \\
          4 & 1 & 8 \\
          -4 & 8 & 1 \\
        \end{array}
      \right]
$$}

 {\small
$$T_{11}:= \frac{1}{11}
\left[ \begin{array}{rrr}
          2 & -9 & -6 \\
          9 & -2 & 6 \\
          -6 & -6 & 7 \\
        \end{array}
      \right], T_{13}:= \frac{1}{13}
\left[ \begin{array}{rrr}
          -4 & -12 & -3 \\
          12 & -3 & -4 \\
          3 & -4 & 12 \\
        \end{array}
      \right],
      \hat{T}_{13}:= \frac{1}{13}
\left[ \begin{array}{rrr}
          0 & -13 & 0 \\
          12 & 0 & 5 \\
          -5 & 0 & 12 \\
        \end{array}
      \right].
$$}

The next matrix can be obtained by multiplying $T_3$ with $T_5$.
We notice a multiplicative structure on this set of matrices. For
the next two prime sizes we have essentially two orthogonal
matrices.

 {\small
$$T_{17}:= \frac{1}{17}
\left[ \begin{array}{rrr}
          12 & -8 & -9 \\
          12 & 9 & 8 \\
          1 & -12 & 12 \\
        \end{array}
      \right], \hat{T_{17}}:= \frac{1}{17}
\left[ \begin{array}{rrr}
          15 & 0 & 8 \\
          8 & 0 & -15 \\
          0 & -17 & 0 \\
        \end{array}
      \right],
$$}

 {\small
$$T_{19}:= \frac{1}{19}
\left[ \begin{array}{rrr}
          6 & -18 & 1 \\
          17 & 6 & 6 \\
          -6 & -1 & 18 \\
        \end{array}
      \right],
      \hat{T}_{19}:= \frac{1}{19}
\left[ \begin{array}{rrr}
          15 & -6 & -10 \\
          10 & 15 & 6 \\
          6 & -10 & 15 \\
        \end{array}
      \right].
$$}

\n From $T_5$, $\hat{T}_{13}$,  and $\hat{T}_{19}$ it is clear
that there is a natural imbedding of the primitive Pythagorean
Triples into this sequence of orthogonal matrices (well known in
the literature). One interesting question here is the following:
what is the algebraic relevance for the geometric invariants
$\alpha_0$, $\alpha$, $\beta$ and $\gamma$, for an orthogonal
matrix as above? One more observation here: each cube in the list
generates by translations and rotations cubes in $C_k$ and two
different cubes cannot generate the same cube because by doing
those transformations the four planes given by the diagonals are
preserved and those are different for two different cubes in the
list. So, in order to count all the cubes in $C_k$, we first
compute the list of irreducible cubes in $\cal L$, up to the side
length $k$, and then for each one we use the formula
(\ref{oldformula}) to find out how many are generated by each in
$C_k$. This is not enough though because there are cubes in $C_k$
which are not irreducible. So, in the end we multiply each cube in
the list $\cal L$ by an integer factor in such a way the resulting
cube can stil fit in $C_k$. Then, we recalculate the invariants on
this cube and use the same formula (\ref{oldformula}) to find the
contribution of the reducible cubes. In the end, we add up all
these numbers and that gives, $NC(k)$, the number of cubes in
$C_k$ .

The first 100 values of $NC$ are: 1, 9, 36, 100, 229, 473, 910,
1648, 2795, 4469, 6818, 10032, 14315, 19907, 27190, 36502, 48233,
62803, 80736, 102550, 128847, 160271, 197516, 241314, 292737,
352591, 421764, 501204, 592257, 696281, 814450, 948112, 1098607,
1267367, 1456292, 1666998, 1901633, 2162179, 2450440, 2768346,
3117935, 3501389, 3923178, 4384792, 4889323, 5439155, 6037660,
6687358, 7391669, 8154671, 8979750, 9870158, 10830095, 11862711,
12972046, 14161848, 15436931, 16801993, 18263634, 19825948,
21493019, 23269647, 25160816, 27171482, 29308957, 31577319,
33986616, 36540004, 39244371, 42106267, 45131996, 48327502,
51700279, 55258019, 59011634, 62965766, 67132037, 71515527,
76127374, 80973598, 86062187, 91401297, 96999986, 102866282,
109014085, 115457359, 122206348, 129266410, 136648555, 144364071,
152426724, 160843660, 169626467, 178787563, 188347314, 198309846,
208694461, 219509943, 230767760, and 242483634.

\section{The code}

The code is written in Maple and it is attached in pdf format  at
the end of the paper after bibliography.

\end{document}


\pagestyle{empty}
\begin{maplegroup}
\begin{_pstyle4}
\begin{mapleinput}
\mapleinline{active}{1d}{restart:with(numtheory):with(plots):}{%
}
\end{mapleinput}
\end{_pstyle4}

\end{maplegroup}

\section{(Step 1) Procedures to find the possible values of k, m and n,
and normal plane vector
(a,b,c)}

\begin{maplegroup}
\begin{_pstyle6}
{\large Find all values of k less than n such that k is of the form a
product of primes of the form 3k+1}
\end{_pstyle6}

\end{maplegroup}
\begin{maplegroup}
\begin{_pstyle4}
\begin{mapleinput}
\mapleinline{active}{1d}{kvalues:=proc(n)
local i,j,k,L,a,p,q,r,m,mm;
L:=\{\};mm:=floor((n+1)/2);
for i from 2 to mm do
a:=ifactors(2*i-1);
k:=nops(a[2]);r:=0;
    for j from 1 to k do
      m:=a[2][j][1]; p:=m mod 3;
      if m=3 then r:=1 fi;
      if r=0 and p=2  then r:=1 fi;
    od;
if r=0 then L:=L union \{2*i-1\};fi;
od;
L:=L union \{1\};
L:=convert(L,list);
end:}{%
}
\end{mapleinput}
\end{_pstyle4}

\end{maplegroup}
\begin{maplegroup}
\begin{_pstyle6}
\QTR{_cstyle9}{Find solutions (m,n)  of the equation
k\symbol{94}2=m\symbol{94}2-mn+n\symbol{94}2,
gcd(m,n)=1,m\TEXTsymbol{>}0,n\TEXTsymbol{>}0,2m\TEXTsymbol{<}n.}
\end{_pstyle6}

\end{maplegroup}
\begin{maplegroup}
\begin{_pstyle4}
\begin{mapleinput}
\mapleinline{active}{1d}{listofmn:=proc(k)
local a,b,i,nx,x,m,n,L,b1,b2;
x:=[isolve(k^2=m^2-m*n+n^2)];
nx:=nops(x);
L:=\{\};
for i from 1 to nx do
if lhs(x[i][1])=m then a:=rhs(x[i][1]); b:=rhs(x[i][2]);
else  b:=rhs(x[i][1]); a:=rhs(x[i][2]); fi;
 if  gcd(a,b)=1 and a>=0 and b>0 and 2*a<b then L:=L union
\{[a,b]\};fi;
od;L;
end:}{%
}
\end{mapleinput}
\end{_pstyle4}

\end{maplegroup}
\begin{maplegroup}
\end{maplegroup}
\begin{maplegroup}
\begin{flushleft}
\QTR{_cstyle9}{Find Only the solutions that satisfy gcd(a,b,c)=1 and
0\TEXTsymbol{<}a\TEXTsymbol{<}=b\TEXTsymbol{<}=c are in the input.}
\end{flushleft}

\end{maplegroup}
\begin{maplegroup}
\begin{mapleinput}
\mapleinline{active}{1d}{abcsol:=proc(d)
local i,j,k,m,u,x,y,sol,cd;sol:=\{\};
for i from 1 to d do
 u:=[isolve(3*d^2-i^2=x^2+y^2)];k:=nops(u);
   for j from 1 to k do
     if rhs(u[j][1])>=i and rhs(u[j][2])>=i then
        cd:=gcd(gcd(i,rhs(u[j][1])),rhs(u[j][2]));
      if cd=1 then sol:=sol union
\{sort([i,rhs(u[j][1]),rhs(u[j][2])])\};fi;
     fi;
   od;
 od;
convert(sol,list);
end:}{%
}
\end{mapleinput}

\end{maplegroup}

\section{(Step2)  The new method of finding r and s. It is followed by
the implementation of the parametrization of equilateral triangles
constructing only one regular tetrahedron and completing it to a cube
(the origin is one of its vertices)}

\begin{maplegroup}
\begin{flushleft}
\textit{{\large Procedure to get the unique writing as
x\symbol{94}2+3y\symbol{94}2 for a prime p of the form 6k+1.}}
\end{flushleft}

\end{maplegroup}
\begin{maplegroup}
\begin{mapleinput}
\mapleinline{active}{1d}{uniquedecomposition:=proc(p)
local out,s,out1,out2;
if p=2 then out:=[1,1]; fi;
 if p>2 then s:=isolve(p=x^2+3*y^2);
  out1:=abs(rhs(s[1][1]));out2:=abs(rhs(s[1][2]));
   if out1^2+3*out2^2=p then out:=[out1,out2]; else
out:=[out2,out1];fi;
 fi;
out; end:}{%
}
\end{mapleinput}

\end{maplegroup}
\begin{maplegroup}
\begin{flushleft}
\textit{{\large The next procedure is calculating the
factorization in Z[sqrt(3)i] of a number of the form u+v*sqrt(3)i
}}
\end{flushleft}

\end{maplegroup}
\begin{maplegroup}
\begin{mapleinput}
\mapleinline{active}{1d}{factoroverEisensteinintegers:=proc(u,v)
local i,N,M,k,a,x,f1,f2,g,y,y1,y2,L,NN,MM,uu,vv;
a:=sqrt(3)*I;
  NN:=gcd(u,v);uu:=u/NN;vv:=v/NN;
  N:=uu^2+3*vv^2;
    x:=uu+vv*a;
  M:=ifactors(N);k:=nops(M[2]);
 for i from 1 to k do
   f1:=M[2][i][1]; f2:=M[2][i][2];
if f1>2 then
   g:=uniquedecomposition(f1);
   else  g:=[1,1]; f2:=1;
fi;
   y:=expand(rationalize(x/(g[1]+a*g[2])));
   y1:=Re(y);
   y2:=type(y1,integer);
     if y2=true then L[i]:=[g[1]+a*g[2],f2]; else
       L[i]:=[g[1]-a*g[2],f2];
     fi;
od;
[NN,seq(L[i],i=1..k),expand(NN*product(L[ii][1]^L[ii][2],ii=1..k))];
end:}{%
}
\end{mapleinput}

\end{maplegroup}
\begin{maplegroup}
\begin{flushleft}
\textit{{\large Finding the gcd between A+isqrt(3)B and 2(q)}}
\end{flushleft}

\end{maplegroup}
\begin{maplegroup}
\begin{mapleinput}
\mapleinline{active}{1d}{findgcd:=proc(A,B,q)
local
i,j,f,common,m,qq,fac,a,L,Important,nfac,f1,ff1,f2,s,LL,n,P,rs,LLL;
a:=sqrt(3)*I;L:=\{\};
f:=factoroverEisensteinintegers(A,B);m:=nops(f)-1;
common:=gcd(f[1],2*q);
qq:=2*q/common^2;
P:=common;
fac:=ifactors(qq);nfac:=nops(fac[2]);
  for i from 1 to nfac do
    f1:=fac[2][i][1];f2:=fac[2][i][2];
       if f1=2 then f2:=1;fi;
         s:=uniquedecomposition(f1);
         for j from 1 to m do
           if s[1]+a*s[2]=f[j][1] then
             P:=P*(s[1]+a*s[2])^(min(f[j][2],f2)); fi;
           if s[1]-a*s[2]=f[j][1] then
             P:=P*(s[1]-a*s[2])^(min(f[j][2],f2)); fi;
                  od;
         od;
P:=expand(P);
rs:=[Re(P),Im(P)/sqrt(3)];
[rs,A*rs[1]+3*B*rs[2] mod 2*q,A*rs[2]-B*rs[1] mod 2*q]
end:}{%
}
\end{mapleinput}

\end{maplegroup}
\begin{maplegroup}
\begin{_pstyle10}
\QTR{_cstyle14}{The fourth point to get a regular tetrahedron is
the origin. We have included only one since all the other are
going to produce the same cube when transformed into the positive
quadrant. }
\end{_pstyle10}

\end{maplegroup}
\begin{maplegroup}
\begin{mapleinput}
\mapleinline{active}{1d}{findpar:=proc(a,b,c,mm,nn)
local i,j,sol,mx,nx,r,s,my,ny,q,d,u,v,w,x,y,z,mu,nu,mv,nv,ef,ns,mz,nz
,mw,nw,T,l,uu,t,R1,R2,fc,k,A,B,rs,G,DD,E,F;

q:=a^2+b^2;k:=sqrt(mm^2-mm*nn+nn^2);d:=sqrt((a^2+b^2+c^2)/3);A:=a*c;B:
=b*d;
   rs:=findgcd(A,-B,q);
   r:=rs[1][2]; s:=rs[1][1];
   uu:=(s^2+3*r^2-2*q)^2;   if uu>0 then t:=s;s:=r;r:=t; fi;
   mx:=-(d*b*(3*r+s)+a*c*(r-s))/(2*q);nx:=-(r*a*c+d*b*s)/q;
   my:=(d*a*(3*r+s)-b*c*(r-s))/(2*q);ny:=-(r*b*c-d*a*s)/q;
   mz:=(r-s)/2;nz:=r;
   mu:=nx;mv:=ny;mw:=nz;nu:=nx-mx;nv:=ny-my; nw:=nz-mz;
   u:=mu*m-nu*n;v:=mv*m-nv*n;w:=mw*m-nw*n;
   x:=mx*m-nx*n;y:=my*m-ny*n;z:=mz*m-nz*n;
   R2:=[(x+u-2*a*k)/3,(v+y-2*b*k)/3,(z+w-2*c*k)/3];
   R1:=[(x+u+2*a*k)/3,(v+y+2*b*k)/3,(z+w+2*c*k)/3];
   fc:=subs(m=mm,n=nn,R1[1]);
if fc=floor(fc) then
 T:=subs(m=mm,n=nn,[[u,v,w],[x,y,z],R1]); else
 T:=subs(m=mm,n=nn,[[u,v,w],[x,y,z],R2]);
fi;
 G:=[(T[1][1]+T[2][1]+T[3][1])/2,(T[1][2]+T[2][2]+T[3][2])/2,(T[1][3]+
T[2][3]+T[3][3])/2];
 DD:=[(T[1][1]+T[2][1]-T[3][1])/2,(T[1][2]+T[2][2]-T[3][2])/2,(T[1][3]
+T[2][3]-T[3][3])/2];
 E:=[(T[1][1]-T[2][1]+T[3][1])/2,(T[1][2]-T[2][2]+T[3][2])/2,(T[1][3]-
T[2][3]+T[3][3])/2];
 F:=[(-T[1][1]+T[2][1]+T[3][1])/2,(-T[1][2]+T[2][2]+T[3][2])/2,(-T[1][
3]+T[2][3]+T[3][3])/2];
[T[1],T[2],T[3],G,DD,E,F];
end:}{%
}
\end{mapleinput}

\end{maplegroup}
\emptyline

\section{(Step 3) Procedures for adding vectors, other operations on
vertices, finding the four k-values of a cube (most of the time is
less than four)}

\begin{maplegroup}
\begin{mapleinput}
\mapleinline{active}{1d}{subtrv:=proc(U,V)
local W;
W[1]:=U[1]-V[1];W[2]:=U[2]-V[2];W[3]:=U[3]-V[3];
[W[1],W[2],W[3]];
end:}{%
}
\end{mapleinput}

\end{maplegroup}
\begin{maplegroup}
\begin{mapleinput}
\mapleinline{active}{1d}{addvec:=proc(U,V)
local W;
W[1]:=U[1]+V[1];W[2]:=U[2]+V[2];W[3]:=U[3]+V[3];
[W[1],W[2],W[3]];
end:}{%
}
\end{mapleinput}

\end{maplegroup}
\begin{maplegroup}
\begin{mapleinput}
\mapleinline{active}{1d}{addvectCube:=proc(T,v)
local i,Q;Q:=\{\};
for i from 1 to 8 do
Q:=Q union \{addvec(T[i],v)\};
od;
Q;
end:}{%
}
\end{mapleinput}

\end{maplegroup}
\begin{maplegroup}
\begin{mapleinput}
\mapleinline{active}{1d}{distance:=proc(A,B)
local C;
C:=subtrv(A,B);
sqrt(C[1]^2+C[2]^2+C[3]^2);
end:}{%
}
\end{mapleinput}

\end{maplegroup}
\begin{maplegroup}
\begin{mapleinput}
\mapleinline{active}{1d}{finddiagonal1:=proc(C)
local i,m,L,x,y,CC;
  L:=convert(C,list);m:=nops(L);
for i from 1 to m do
x:=type(distance(L[1],L[i])/sqrt(3),integer);
if x=true then y:=[L[1],L[i]];fi;
od;
CC:=C minus \{y[1],y[2]\};
y,CC;
end:}{%
}
\end{mapleinput}

\end{maplegroup}
\begin{maplegroup}
\begin{mapleinput}
\mapleinline{active}{1d}{finddiagonal:=proc(C)
local x,x1,x2;
x:=finddiagonal1(C);
x1:=finddiagonal1(x[2]);
x2:=finddiagonal1(x1[2]);
[x[1],x1[1],x2[1],convert(x2[2],list)];
end:}{%
}
\end{mapleinput}

\end{maplegroup}
\begin{maplegroup}
\begin{mapleinput}
\mapleinline{active}{1d}{unitvector:=proc(U)
local i,j,k,l,x;
i:=U[1];j:=U[2];k:=U[3];
l:=gcd(gcd(i,j),k);
x:=(i^2+j^2+k^2)/(3*l^2);
sqrt(x);
end:}{%
}
\end{mapleinput}

\end{maplegroup}
\begin{maplegroup}
\begin{_pstyle4}
\begin{mapleinput}
\mapleinline{active}{1d}{fourkvalues:=proc(C)
local N1,N2,N3,N4,x,length;length:=distance(C[1],C[2]);
x:=finddiagonal(C);
N1:=unitvector(subtrv(x[1][1],x[1][2]));
N2:=unitvector(subtrv(x[2][1],x[2][2]));
N3:=unitvector(subtrv(x[3][1],x[3][2]));
N4:=unitvector(subtrv(x[4][1],x[4][2]));
\{length/N1,length/N2,length/N3,length/N4\};
end:}{%
}
\end{mapleinput}
\end{_pstyle4}

\end{maplegroup}
\begin{maplegroup}
\begin{mapleinput}
\mapleinline{active}{1d}{multbyfactorv:=proc(v,k)
local w;
w[1]:=v[1]*k;w[2]:=v[2]*k;w[3]:=v[3]*k;
[w[1],w[2],w[3]];
end:}{%
}
\end{mapleinput}

\end{maplegroup}
\begin{maplegroup}
\begin{_pstyle4}
\begin{mapleinput}
\mapleinline{active}{1d}{multbyfactor:=proc(T,k)
local i,NT,Q;NT:=\{\};
Q:=convert(T,list);
for i from 1 to 8 do
NT:=NT union \{multbyfactorv(Q[i],k)\};
od;NT;
end:}{%
}
\end{mapleinput}
\end{_pstyle4}

\end{maplegroup}

\section{(Step 4) Translation of the cube obtained in Step 2  into the
first quadrant and getting the orbit and the generalized orbit, then
intersections with its translations along (0,0,1) and (0,1,-1).
Calculation of alpha\_0, alpha, beta and gamma for a cube.}

\begin{maplegroup}
\begin{flushleft}
\textit{{\large Translating the cube to the first octant
(x\TEXTsymbol{>}=0,y\TEXTsymbol{>}=0,z\TEXTsymbol{>}=0)}}
\end{flushleft}

\end{maplegroup}
\begin{maplegroup}
\begin{mapleinput}
\mapleinline{active}{1d}{tmttopqcube:=proc(T)
local i,a,b,c,v,O,C;
a:=min(T[1][1],T[2][1],T[3][1],T[4][1],T[5][1],T[6][1],T[7][1],0);
b:=min(T[1][2],T[2][2],T[3][2],T[4][2],T[5][2],T[6][2],T[7][2],0);
c:=min(T[1][3],T[2][3],T[3][3],T[4][3],T[5][3],T[6][3],T[7][3],0);
O:=[0,0,0];
v:=[a,b,c];C:=
\{subtrv(O,v),subtrv(T[1],v),subtrv(T[2],v),subtrv(T[3],v),
subtrv(T[4],v),subtrv(T[5],v),subtrv(T[6],v),subtrv(T[7],v)\};
end:}{%
}
\end{mapleinput}

\end{maplegroup}
\begin{maplegroup}
\begin{mapleinput}
\mapleinline{active}{1d}{mscofmcube:=proc(Q)
local a,b,c,T;
T:=convert(Q,list);
a:=max(T[1][1],T[2][1],T[3][1],T[4][1],T[5][1],T[6][1],T[7][1],T[8][1]
);
b:=max(T[1][2],T[2][2],T[3][2],T[4][2],T[5][2],T[6][2],T[7][2],T[8][2]
);
c:=max(T[1][3],T[2][3],T[3][3],T[4][3],T[5][3],T[6][3],T[7][3],T[8][3]
);
max(a,b,c);
end:}{%
}
\end{mapleinput}

\end{maplegroup}
\begin{maplegroup}
\begin{mapleinput}
\mapleinline{active}{1d}{orbit1Cube:=proc(T)
local
i,k,T1,a,b,c,x,T2,T3,T4,T5,T6,T7,T8,T9,T10,T11,T12,T13,T14,T15,T16,T17
,T18,
T19,T20,T21,T22,T23,T24,S,Q,d,a1,b1,c1;
Q:=convert(T,list);
d:=mscofmcube(T);
T1:=T;
T2:=\{seq([Q[k][2],Q[k][3],Q[k][1]],k=1..8)\};
T3:=\{seq([Q[k][1],Q[k][3],Q[k][2]],k=1..8)\};
T4:=\{seq([Q[k][1],Q[k][2],d-Q[k][3]],k=1..8)\};
T5:=\{seq([Q[k][2],Q[k][3],d-Q[k][1]],k=1..8)\};
T6:=\{seq([Q[k][1],Q[k][3],d-Q[k][2]],k=1..8)\};
T7:=\{seq([Q[k][1],d-Q[k][2],Q[k][3]],k=1..8)\};
T8:=\{seq([Q[k][2],d-Q[k][3],Q[k][1]],k=1..8)\};
T9:=\{seq([Q[k][1],d-Q[k][3],Q[k][2]],k=1..8)\};
T10:=\{seq([d-Q[k][1],Q[k][2],Q[k][3]],k=1..8)\};
T11:=\{seq([d-Q[k][2],Q[k][3],Q[k][1]],k=1..8)\};
T12:=\{seq([d-Q[k][1],Q[k][3],Q[k][2]],k=1..8)\};
T13:=\{seq([Q[k][1],d-Q[k][2],d-Q[k][3]],k=1..8)\};
T14:=\{seq([Q[k][2],d-Q[k][3],d-Q[k][1]],k=1..8)\};
T15:=\{seq([Q[k][1],d-Q[k][3],d-Q[k][2]],k=1..8)\};
T16:=\{seq([d-Q[k][1],d-Q[k][2],Q[k][3]],k=1..8)\};
T17:=\{seq([d-Q[k][2],d-Q[k][3],Q[k][1]],k=1..8)\};
T18:=\{seq([d-Q[k][1],d-Q[k][3],Q[k][2]],k=1..8)\};
T19:=\{seq([d-Q[k][1],Q[k][2],d-Q[k][3]],k=1..8)\};
T20:=\{seq([d-Q[k][2],Q[k][3],d-Q[k][1]],k=1..8)\};
T21:=\{seq([d-Q[k][1],Q[k][3],d-Q[k][2]],k=1..8)\};
T22:=\{seq([d-Q[k][1],d-Q[k][2],d-Q[k][3]],k=1..8)\};
T23:=\{seq([d-Q[k][2],d-Q[k][3],d-Q[k][1]],k=1..8)\};
T24:=\{seq([d-Q[k][1],d-Q[k][3],d-Q[k][2]],k=1..8)\};
S:=
\{T1,T2,T3,T4,T5,T6,T7,T8,T9,T10,T11,T12,T13,T14,T15,T16,T17,T18,T19,T
20,T21,T22,T23,T24\};
S;
end:}{%
}
\end{mapleinput}

\end{maplegroup}
\begin{maplegroup}
\begin{mapleinput}
\mapleinline{active}{1d}{orbitCube:=proc(T)
local S,Q,T1;
Q:=convert(T,list);
T1:=\{seq([Q[k][3],Q[k][2],Q[k][1]],k=1..8)\};
S:=orbit1Cube(T) union orbit1Cube(T1);
S;
end:}{%
}
\end{mapleinput}

\end{maplegroup}
\begin{maplegroup}
\begin{mapleinput}
\mapleinline{active}{1d}{translCube:=proc(T)
local S,Q,i,j,k,a2,b2,c2,a,b,c,d,alpha;
Q:=convert(T,list);
a:=max(seq(Q[k][1],k=1..8));
b:=max(seq(Q[k][2],k=1..8));
c:=max(seq(Q[k][3],k=1..8));
d:=max(a,b,c);
 a2:=d-a;b2:=d-b;c2:=d-c;
  S:=orbitCube(T);alpha:=nops(S);
 for i from 0 to a2 do
   for j from 0 to b2 do
    for k from 0 to c2 do
        S:=S union orbitCube(addvectCube(T,[i,j,k]));
      od;
   od;
 od;
#alpha,nops(S),[a2,b2,c2],(a2*b2+b2*c2+c2*a2)*alpha;
S;
end:}{%
}
\end{mapleinput}

\end{maplegroup}
\begin{maplegroup}
\begin{mapleinput}
\mapleinline{active}{1d}{alphabetagamma:=proc(T)
local x,S,i,S1,S2,d,cube,j,exitv,alpha,beta,gamma;
d:=mscofmcube(T);x:=nops(orbitCube(T));
S1:=translCube(T);S2:=S1;S:=convert(S2,list);alpha:=nops(S);
 for i from 1 to alpha do
  cube:=S[i];
  j:=1;  exitv:=0;
   while j<=8 and exitv=0 do
     if cube[j][3]=d then exitv:=1;
     S2:=S2 minus \{cube\};fi;
     j:=j+1;
   od;
 od;
beta:=nops(S2);S:=convert(S2,list);
for i from 1 to beta do
  cube:=S[i];
  j:=1;  exitv:=0;
   while j<=8 and exitv=0 do
     if cube[j][2]=d then exitv:=1;
     S2:=S2 minus \{cube\};fi;
     j:=j+1;
   od;
 od;
gamma:=nops(S2);
[S2,[x,alpha,beta,gamma]];
end:}{%
}
\end{mapleinput}

\end{maplegroup}
\emptyline

\section{\QTR{_cstyle13}{(Step 5) The function that gives the number of cubes generated by a
small cube which fits in a cube of dimension d,  in an arbitrary
bigger cube of dimension n.}}

\begin{maplegroup}
\begin{_pstyle4}
\begin{mapleinput}
\mapleinline{active}{1d}{f:=(n,d,alpha,beta,gamma)->(n-d+1)^3*alpha-3*(n-d)*(-d+1+n)^2*beta+3*
gamma*(n-d+1)*(n-d)^2:}{%
}
\end{mapleinput}
\end{_pstyle4}

\end{maplegroup}

\section{\QTR{_cstyle13}{(Step 6) List of minimal (irreducible) cubes together with n, minimal
size containg the cube, the k-values, alpha, beta, and gamma. In the
end another list is created in a new procedure to get the list of the
reducible cubes together with  the alpha, beta, and gamma.}}

\begin{maplegroup}
\begin{_pstyle10}
{\large We compute the list of minimal cubes taking into account that
at some point only k=1 is not sufficient to generate a cube . }
\end{_pstyle10}

\end{maplegroup}
\begin{maplegroup}
\begin{_pstyle4}
\begin{mapleinput}
\mapleinline{active}{1d}{ExtendList:=proc(n,N,L::list,mm,nn,Orb::array)
local i,sol,nsol,nel,ttpcube,orb::array,
LL::list,tnel,NL::list,cio,C,pC,k,kvalues,m,exception,abg,NN;
NN:=floor((N-1)/2);
orb:=array(1..2*NN+1);
nel:=nops(L);
LL:=L;
k:=sqrt(mm^2-mm*nn+nn^2);
m:=n*k;
if m<=N then
 sol:=abcsol(n);nsol:=nops(sol);
 tnel:=nel;orb:=Orb;
   for i from 1 to nsol do

C:=findpar(sol[i][1],sol[i][2],sol[i][3],mm,nn);pC:=tmttopqcube(C);

        kvalues:=fourkvalues(pC);exception:=evalb(1 in kvalues);

           if k=1 or exception=false then
           cio:=evalb(pC in orb[m]);
           if cio=false then
              abg:=alphabetagamma(pC);
              orb[m]:=orb[m] union abg[1];

NL[tnel+1]:=[n,mscofmcube(pC),pC,kvalues,abg[2],sol[i]];
             tnel:=tnel+1;
           fi;
      fi;
      od;
LL:=[seq(L[i],i=1..nel),seq(NL[j],j=nel+1..tnel)];
fi;
LL,orb;
end:
}{%
}
\end{mapleinput}
\end{_pstyle4}

\end{maplegroup}
\begin{maplegroup}
\begin{mapleinput}
\mapleinline{active}{1d}{ExtendListuptoN:=proc(N)
local i,j,k,l,kv,kvn,NN,L,Orb::array,mn,nmn,E,n,ii;
 kv:=kvalues(N);kvn:=nops(kv);
  NN:=floor((N-1)/2);
  L:=[];Orb:=array(1..2*NN+1);
      for ii from 1 to 2*NN+1 do
        Orb[ii]:=\{\};
      od;
    for i from 1 to kvn do
     k:=kv[i];mn:=listofmn(k);nmn:=nops(mn);
      for j from 1 to nmn do
       for l from 1 to NN+1 do
         n:=2*l-1;
         E:=ExtendList(n,N,L,mn[j][1],mn[j][2],Orb);
              L:=E[1];
        for ii from 1 to 2*NN+1 do
           Orb[ii]:=E[2][ii];
        od;

   od;
  od;
 od;
L;
end:}{%
}
\end{mapleinput}

\end{maplegroup}
\begin{maplegroup}
\begin{mapleinput}
\mapleinline{active}{1d}{L:=ExtendListuptoN(50):}{%
}
\end{mapleinput}

\end{maplegroup}
\begin{maplegroup}
\begin{mapleinput}
\mapleinline{active}{1d}{ExtendListuptoNmultiples:=proc(N,L)
local i,j,x,lc,m,mm,d,dd,C,CC,LL;
m:=nops(L);i:=1;LL:=\{\};
while i<=m  do
  d:=L[i][2];
     if d<=N then
        mm:=floor(N/d);C:=L[i][3];j:=2;
         while j<=mm  do
             CC:=multbyfactor(C,j);x:=alphabetagamma(CC);
              dd:=d*j;lc:=nops(LL);
              LL:=LL union \{[L[i][1]*j,dd,CC,L[i][4],x[2]]\};
             j:=j+1;
         od;
      fi;
    i:=i+1;
 od;
convert(LL,list);
end:}{%
}
\end{mapleinput}

\end{maplegroup}
\begin{maplegroup}
\begin{mapleinput}
\mapleinline{active}{1d}{LL:=ExtendListuptoNmultiples(50,L):}{%
}
\end{mapleinput}

\end{maplegroup}
\begin{maplegroup}
\end{maplegroup}

\section{(Step 7) Adding up the contribution of each cube.}

\begin{maplegroup}
\end{maplegroup}
\begin{maplegroup}
\begin{mapleinput}
\mapleinline{active}{1d}{addupnew:=proc(N,L,LL)
local i,j,k,nc,mm,m,d,C,CC,x,dd,nt;
nc:=0;
m:=nops(L);
i:=1;
  while i<=m  do
  d:=L[i][2];
   if d<=N then
             nc:=nc+f(N,d,L[i][5][2],L[i][5][3],L[i][5][4]);
   fi;  i:=i+1;
  od;
m:=nops(LL);
i:=1;
  while i<=m  do
  d:=LL[i][2];
   if d<=N then
             nc:=nc+f(N,d,LL[i][5][2],LL[i][5][3],LL[i][5][4]);
   fi;  i:=i+1;
  od;
nc;
end:
}{%
}
\end{mapleinput}

\end{maplegroup}
\begin{maplegroup}
\end{maplegroup}
\begin{maplegroup}
\begin{mapleinput}
\mapleinline{active}{1d}{NC:=[seq(addupnew(k,L,LL),k=1..50)];}{%
}
\end{mapleinput}

\mapleresult
\begin{maplelatex}
\mapleinline{inert}{2d}{NC := [1, 9, 36, 100, 229, 473, 910, 1648, 2795, 4469, 6818, 10032,
14315, 19907, 27298, 36886, 49133, 64531, 83784, 107542, 136551,
171599, 213524, 263202, 321849, 390415, 469932, 561492, 667305,
789317, 929098, 1088500, 1269367, 1473635, 1703708, 1961706, 2251289,
2575291, 2936272, 3337026, 3780455, 4269605, 4813854, 5414560,
6076915, 6804587, 7603120, 8476390, 9430481, 10471175];}{%
\maplemultiline{
\mathit{NC} := [1, \,9, \,36, \,100, \,229, \,473, \,910, \,1648
, \,2795, \,4469, \,6818, \,10032, \,14315, \,19907,  \\
27298, \,36886, \,49133, \,64531, \,83784, \,107542, \,136551, \,
171599, \,213524, \,263202,  \\
321849, \,390415, \,469932, \,561492, \,667305, \,789317, \,
929098, \,1088500, \,1269367,  \\
1473635, \,1703708, \,1961706, \,2251289, \,2575291, \,2936272,
\,3337026, \,3780455,  \\
4269605, \,4813854, \,5414560, \,6076915, \,6804587, \,7603120,
\,8476390, \,9430481,  \\
10471175] }
}
\end{maplelatex}

\end{maplegroup}
\begin{maplegroup}
\begin{mapleinput}
\end{mapleinput}

\end{maplegroup}
\begin{maplegroup}
\begin{mapleinput}
\end{mapleinput}

\end{maplegroup}
\begin{maplegroup}
\end{maplegroup}
\begin{maplegroup}
\end{maplegroup}
\begin{maplegroup}
\end{maplegroup}
\begin{maplegroup}
\end{maplegroup}
\begin{maplegroup}
\end{maplegroup}
\begin{maplegroup}
\end{maplegroup}
\begin{maplegroup}
\end{maplegroup}
\begin{maplegroup}
\end{maplegroup}
\begin{maplegroup}
\end{maplegroup}
\begin{maplegroup}
\end{maplegroup}
\begin{maplegroup}
\end{maplegroup}
\begin{maplegroup}
\end{maplegroup}
\begin{maplegroup}
\end{maplegroup}
\begin{maplegroup}
\end{maplegroup}
\begin{maplegroup}
\end{maplegroup}
\begin{maplegroup}
\end{maplegroup}
\begin{maplegroup}
\end{maplegroup}
\begin{maplegroup}
\end{maplegroup}
\emptyline